\documentclass{tran-l}

\newtheorem{thm}{Theorem}[subsection]
\newtheorem{cor}[thm]{Corollary}
\newtheorem{lem}[thm]{Lemma}

\theoremstyle{definition}

\theoremstyle{remark}

\numberwithin{equation}{subsection}

\begin{document}

\title[Binary matrices associated with the Tribonacci sequence.]
 {Inverse Powers and $2$-Adic Properties of Binary $3\times 3$ Matrices 
Associated with the Tribonacci Sequence.}

\author{Wilson Arley Martinez * , Samin Ingrith Ceron}

\address{Martinez, W.A.; Departmento de Matem\'aticas, Universidad del Cauca , Popay\'an , Colombia}
\email{wamartinez@unicauca.edu.co}

\address{Ceron, S.I.; Departmento de Matem\'aticas, Universidad del Cauca , Popay\'an , Colombia}
\email{sicbravo@gmail.com}

\thanks{This work was completed with the support of the Universidad del Cauca.}

\thanks{The author was also supported by the research group “Estructuras Algebraicas, Divulgación Matemática y Teorías Asociadas. @DiTa”.}

\thanks{* Corresponding Author: Wilson Arley Martinez.}




\date{May 25, 2026 and, in revised form, May, 2026.}

\dedicatory{}

\commby{W.A.M}


\begin{abstract}
In this paper, we study binary $3\times 3$ matrices with determinant one that generate the Tribonacci recurrence. Starting from a specific generating matrix, we obtain explicit expressions for its powers and inverse powers in terms of the Tribonacci sequence $T_n$ and two associated sequences $B_n$ and $R_n$.  This matrix framework yields identities involving traces, determinants, and products, including addition formulas and Cassini-type relations linking $T_n$, $B_n$, and $R_n$. Moreover, we derive new identities from the matrix representation of inverse powers and establish structural relations between these sequences. As an application, we investigate arithmetic properties of Tribonacci numbers and determine the exact $2$-adic valuation of $T_{2^k}$ for sufficiently large $k$. Finally, we provide a classification of $3\times 3$ matrices associated with the Tribonacci recurrence via similarity over $\mathbb{Q}$.
\end{abstract}

\maketitle

\section*{Introduction}

\section{Introduction}

Linear recurrence sequences play a central role in number theory, combinatorics, and linear algebra. Among them, the \emph{Tribonacci sequence} $T_n$ is a natural generalization of the Fibonacci sequence and is defined by the third--order linear recurrence relation
\begin{equation*}
    T_{n+3} = T_{n+2} + T_{n+1} + T_n, \quad \text{for } n \ge 0,
\end{equation*}
with initial values
\[
T_0 = 0, \quad T_1 = 0, \quad T_2 = 1.
\]
The first terms of the sequence are
\[
0,\; 0,\; 1,\; 1,\; 2,\; 4,\; 7,\; 13,\; 24,\; 44,\; \dots
\]
This sequence has been widely studied due to its rich algebraic and arithmetic structure, including closed-form expressions, divisibility properties, and periodic behavior modulo integers.

Matrix methods provide a powerful framework for studying higher--order linear recurrence relations. It is well known (see, for example, \cite{GravesBaconCook2013}) that the Tribonacci sequence can be generated using powers of the companion matrix
\[
A =
\begin{bmatrix}
1 & 1 & 1 \\[1mm]
1 & 0 & 0 \\[1mm]
0 & 1 & 0
\end{bmatrix}.
\]
More precisely, for all integers $n \ge 2$, one has
\[
A^n =
\begin{bmatrix}
T_{n+2} & T_{n+1}+ T_n     & T_{n+1} \\[1mm]
T_{n+1} & T_{n}+ T_{n-1}   & T_{n} \\[1mm]
T_n     & T_{n-1}+ T_{n-2} & T_{n-1}
\end{bmatrix},
\]
where the indices are interpreted consistently with the initial conditions. This representation allows the Tribonacci sequence to be studied through algebraic properties of matrices, such as determinants, traces, and conjugation.

\quad

While the classical approach is well developed (see, for example, \cite{Feng2011, ShannonHoradam1972,  Spickerman1982}), some structural aspects of matrix representations of the Tribonacci sequence remain unexplored, especially for binary matrices and inverse powers. In this paper, we study binary \(3\times 3\) matrices with determinant one that generate the Tribonacci recurrence. From a canonical matrix, we derive explicit formulas for both positive and negative powers.

\quad

These representations involve the Tribonacci sequence \(T_n\) and two related sequences: a Tribonacci-type sequence \(B_n\) and a signed version \(R_n\), providing a unified framework for identities and structural analysis. Using this approach, we obtain addition formulas, trace and determinant relations, and Cassini-type identities connecting \(T_n\), \(B_n\), and \(R_n\). Inverse powers yield new nontrivial identities involving \(R_n\).

\quad

As an application, we study arithmetic properties of \(T_n\), determining the exact \(2\)-adic valuation of \(T_{2^k}\) for all sufficiently large \(k\), refining previous results (see, for example, \cite{MarquesLengyel2014,Sobolewski2017}). Finally, we classify matrices associated with the recurrence, showing that every \(3\times 3\) matrix over \(\mathbb{Q}\) with characteristic polynomial \(x^3 - x^2 - x - 1\) is similar to a companion matrix; hence all such matrices are conjugate in \(GL(3,\mathbb{Q})\).

\quad

The paper is organized as follows: Section~2 develops the matrix framework and identities; Section~3 studies divisibility and the \(2\)-adic valuation of \(T_{2^k}\); the final section presents the classification of the associated matrices.

\section{Recurrences Generated by Binary $3 \times 3$ Matrices with Determinant One}
Consider the generating matrix
\[
u = \begin{pmatrix}0 & 0 & 1\\ 1 & 0 & 1\\ 0 & 1 & 1\end{pmatrix}.
\]
In this section, we present the Tribonacci \textit{u}-matrix with determinant~1, which provides a matrix representation of Tribonacci numbers. We use it to compute powers \(u^n\), determinants,, traces,  inverses, and Cassini-like identities; to derive the characteristic roots and a Binet-type formula for the sequence \( T_n \); to establish identities involving the Tribonacci sequence together with the associated Tribonacci-type sequence \( B_n \) and the signed Tribonacci-type sequence \(R_n\); and to obtain summation formulas for Tribonacci numbers, \( B_n \), and \( R_n \).

\subsection{The Matrix Representation of Recurrences and Their Identities}

\begin{lem}\label{l:3} 
Let \( T_n \) be the Tribonacci numbers defined by
\[
T_n = T_{n-1} + T_{n-2} + T_{n-3}, \quad T_0 = 0,\ T_1 = 0,\ T_2 = 1,
\]
and let \( B_n \) be the sequence defined by
\[
B_n = B_{n-1} + B_{n-2} + B_{n-3}, \quad B_0 = 1,\ B_1 = 0,\ B_2 = 1.
\]
Then, for all \( n \geq 2 \),
\begin{align}
B_{n} &= T_n + T_{n-1},\label{id:1}\\
T_{n+2} &= B_{n}+T_{n+1},\label{id:2}\\
B_{n+1}+B_{n} &= T_{n}+T_{n+2}.\label{id:14}
\end{align}
\end{lem}

\begin{proof}
We prove \eqref{id:1} by strong induction on \( n \ge 2 \).

For the base cases,
\[
B_2 = 1 = T_2 + T_1, \qquad
B_3 = 2 = T_3 + T_2,
\]
so the identity holds.

Assume that \( B_j = T_j + T_{j-1} \) for every integer \( j \) with \( 2 \le j \le k \). Using the recurrence for \(B_n\), we obtain
\begin{align*}
B_{k+1}
    &= B_k + B_{k-1} + B_{k-2} \\
    &= (T_k + T_{k-1}) + (T_{k-1} + T_{k-2}) + (T_{k-2} + T_{k-3}) \\
    &= (T_k + T_{k-1} + T_{k-2}) + (T_{k-1} + T_{k-2} + T_{k-3}) \\
    &= T_{k+1} + T_k,
\end{align*}
where we used the recurrence defining \(T_{k+1}\) and \(T_k\).  
Thus, the identity holds for \( k+1 \), and by strong induction it follows for all \( n \ge 2 \).

\medskip

Identities \eqref{id:2} and \eqref{id:14} follow directly from \eqref{id:1} and the Tribonacci recurrence. Indeed,
\begin{align*}
T_{n+2}
    &= T_{n+1} + T_n + T_{n-1} \\
    &= T_{n+1} + (T_n + T_{n-1}) \\
    &= T_{n+1} + B_n,
\end{align*}
which proves \eqref{id:2}. Moreover,
\begin{align*}
B_{n+1} + B_n
    &= (T_{n+1} + T_n) + (T_n + T_{n-1}) \\
    &= T_n + (T_{n+1} + T_n + T_{n-1}) \\
    &= T_n + T_{n+2},
\end{align*}
establishing \eqref{id:14}.
\end{proof}

\begin{lem}\label{l:1} 
Let \(u\) be the invertible matrix
\[
u = \begin{pmatrix} 
0 & 0 & 1 \\
1 & 0 & 1 \\
0 & 1 & 1
\end{pmatrix}.
\]
Then, for every \( n \ge 1 \),
\[
u^{n} =
\begin{pmatrix}
T_{n-1} & T_{n} & T_{n+1} \\
B_{n-1} & B_{n} & T_{n}+T_{n+1} \\
T_{n} & T_{n+1} & B_{n}+T_{n+1}
\end{pmatrix},
\]
where \(T_n\) is the Tribonacci sequence defined by
\[
T_n = T_{n-1}+T_{n-2}+T_{n-3}, \qquad
T_0=T_1=0,\; T_2=1,
\]
and \(B_n\) satisfies
\[
B_n = B_{n-1}+B_{n-2}+B_{n-3}, \qquad
B_0=1,\; B_1=0,\; B_2=1.
\]
\end{lem}

\begin{proof}
Observe that \(u\) is the companion matrix associated with the Tribonacci recurrence. Hence, one expects the entries of \(u^n\) to be expressed in terms of Tribonacci-type sequences. We verify this claim by induction on \(n\).

For \(n=1\), a direct comparison shows that
\[
u=
\begin{pmatrix}
T_0 & T_1 & T_2 \\
B_0 & B_1 & T_1+T_2 \\
T_1 & T_2 & B_1+T_2
\end{pmatrix},
\]
so the formula holds. Assume the identity is valid for some \(k\ge1\). Then
\[
u^{k+1}=u^k u,
\]
and using the induction hypothesis we compute
\[
u^{k+1}=
\begin{pmatrix}
T_{k-1} & T_k & T_{k+1} \\
B_{k-1} & B_k & T_k+T_{k+1} \\
T_k & T_{k+1} & B_k+T_{k+1}
\end{pmatrix}
\begin{pmatrix}
0 & 0 & 1 \\
1 & 0 & 1 \\
0 & 1 & 1
\end{pmatrix}.
\]
Performing the multiplication gives
\[
u^{k+1}=
\begin{pmatrix}
T_k & T_{k+1} & T_{k-1}+T_k+T_{k+1} \\
B_k & T_k+T_{k+1} & (B_{k-1}+T_k)+(B_k+T_{k+1}) \\
T_{k+1} & B_k+T_{k+1} & (T_k+T_{k+1})+(B_k+T_{k+1})
\end{pmatrix}.
\]
By the Tribonacci recurrence,
\[
T_{k+2}=T_{k-1}+T_k+T_{k+1},
\]
and by Lemma~\ref{l:3}, \(B_n=T_n+T_{n-1}\), which implies
\[
B_{k+1}=T_{k+1}+T_k.
\]
Substituting these identities into the previous matrix yields
\[
u^{k+1}=
\begin{pmatrix}
T_k & T_{k+1} & T_{k+2} \\
B_k & B_{k+1} & T_{k+1}+T_{k+2} \\
T_{k+1} & T_{k+2} & B_{k+1}+T_{k+2}
\end{pmatrix}.
\]
Therefore, the formula holds for \(k+1\), and the result follows by induction.
\end{proof}

Analogously, the following lemma can be proved by mathematical induction.

\begin{lem}\label{l:2} 
Let \( u \) be the invertible matrix  
\[
u=\begin{pmatrix} 
0 & 0 & 1 \\
1 & 0 & 0 \\
1 & 1 & 1
\end{pmatrix}.
\]
Then, for every \( n\in\mathbb{Z}^{+} \),
\[
u^{n} =
\begin{pmatrix}
T_{n}+T_{n-1} & T_{n} & T_{n}+B_{n-1} \\
B_{n-1} & T_{n-1} & T_{n} \\
2T_{n}+B_{n-1} & T_{n}+B_{n-1} & 2T_{n}+T_{n-1}+B_{n-1}
\end{pmatrix},
\]
where \( T_n \) is the Tribonacci sequence defined by
\[
T_n=T_{n-1}+T_{n-2}+T_{n-3},\qquad
T_0=T_1=0,\; T_2=1,
\]
and \(B_n\) is the sequence satisfying
\[
B_n=B_{n-1}+B_{n-2}+B_{n-3},\qquad
B_0=1,\; B_1=0,\; B_2=1.
\]
\end{lem}

\begin{proof}
We prove the result using the principle of mathematical induction.

For \( n=2 \),
\[
u^{2}=
\begin{bmatrix}
T_{2}+T_{1} & T_{2} & T_{2}+B_{1}\\
B_{1} & T_{1} & T_{2}\\
2T_{2}+B_{1} & T_{2}+B_{1} & 2T_{2}+T_{1}+B_{1}
\end{bmatrix}
=
\begin{bmatrix}
1 & 1 & 1\\
0 & 0 & 1\\
2 & 1 & 2
\end{bmatrix},
\]
so the formula holds. Assume that it is true for some \( n=k \), that is,
\[
u^{k}=
\begin{pmatrix}
T_{k}+T_{k-1} & T_{k} & T_{k}+B_{k-1}\\
B_{k-1} & T_{k-1} & T_{k}\\
2T_{k}+B_{k-1} & T_{k}+B_{k-1} & 2T_{k}+T_{k-1}+B_{k-1}
\end{pmatrix}.
\]
Then,
\[
u^{k+1}=u^{k}u
\]
and a direct matrix multiplication gives
\[
\begin{pmatrix}
2T_{k}+B_{k-1} & T_{k}+B_{k-1} & 2T_{k}+T_{k-1}+B_{k-1}\\
T_{k-1}+T_{k} & T_{k} & T_{k}+B_{k-1}\\
3T_{k}+T_{k-1}+2B_{k-1} & 2T_{k}+T_{k-1}+B_{k-1} & 4T_{k}+T_{k-1}+2B_{k-1}
\end{pmatrix}.
\]
Using the recurrence relations for \(T_n\) and \(B_n\) (see Lemma~\ref{l:3}), this matrix becomes
\[
\begin{pmatrix}
T_{k+1}+T_k & T_{k+1} & T_{k+1}+B_k\\
B_k & T_k & T_{k+1}\\
2T_{k+1}+B_k & T_{k+1}+B_k & 2T_{k+1}+T_k+B_k
\end{pmatrix},
\]
which completes the induction.
\end{proof}

\begin{cor}
Let 
\[
u=
\begin{pmatrix} 
0 & 0 & 1 \\
1 & 0 & 1 \\
0 & 1 & 1
\end{pmatrix},
\qquad
v=
\begin{pmatrix} 
0 & 0 & 1 \\
1 & 0 & 0 \\
1 & 1 & 1
\end{pmatrix}.
\]
Both matrices are companion matrices of the Tribonacci polynomial \(x^3-x^2-x-1\); hence they are similar and have the same determinant equal to one. Therefore, for every \( n\in\mathbb{N} \),
\[
\operatorname{Tr}(u^n)=\operatorname{Tr}(v^n)=3B_n+B_{n-1}.
\]
\end{cor}

Although identities of the form
\[
T_{m+n} = T_{m-1}T_{n} + T_{m}(T_{n} + T_{n-1}) + T_{m+1}T_{n+1}
\]
are known (see \cite{Pethe1988}), 
our approach provides a unified derivation via matrix representations.

\begin{cor}\label{c:addition}
For \(m,n\ge1\), the following identities hold:
\begin{align}
T_{m+n} &=  T_{m-1}T_{n} + T_{m}B_{n} + T_{m+1}T_{n+1},\label{id:add1}\\
B_{m+n} &=  B_{m-1}T_{n} + B_{m}B_{n} + (T_{m}+T_{m+1})T_{n+1}. \label{id:add2}  
\end{align}
\end{cor}

\begin{proof}
For \(m,n\ge1\), matrix powers satisfy \(u^{m+n}=u^m u^n\).  
Using the expression for \(u^k\) obtained in Lemma~\ref{l:1}, we write

\[
u^m =
\begin{pmatrix}
T_{m-1} & T_{m} & T_{m+1} \\
B_{m-1} & B_{m} & T_{m}+T_{m+1} \\
T_{m} & T_{m+1} & B_{m}+T_{m+1}
\end{pmatrix},
\quad
u^n =
\begin{pmatrix}
T_{n-1} & T_{n} & T_{n+1} \\
B_{n-1} & B_{n} & T_{n}+T_{n+1} \\
T_{n} & T_{n+1} & B_{n}+T_{n+1}
\end{pmatrix}.
\]
Performing the matrix multiplication gives
\[
u^m u^n =
\begin{pmatrix}
D & A & B \\
E & C & A+B \\
A & B & C+B
\end{pmatrix},
\]
where
\[
\begin{aligned}
A &= T_{m-1}T_{n} + T_{m}B_{n} + T_{m+1}T_{n+1},\\
B &= T_{m-1}T_{n+1} + T_m(T_n+T_{n+1}) + T_{m+1}(T_{n+1}+B_n),\\
C &= B_{m-1}T_{n} + B_{m}B_{n} + (T_{m}+T_{m+1})T_{n+1},\\
D &= T_{m-1}T_{n-1} + T_{m}B_{n-1} + T_{m+1}T_{n},\\
E &= B_{m-1}T_{n-1} + B_{m}B_{n-1} + (T_{m}+T_{m+1})T_{n}.
\end{aligned}
\]
On the other hand,
\[
u^{m+n}=
\begin{pmatrix}
T_{m+n-1} & T_{m+n} & T_{m+n+1} \\
B_{m+n-1} & B_{m+n} & T_{m+n}+T_{m+n+1} \\
T_{m+n} & T_{m+n+1} & B_{m+n}+T_{m+n+1}
\end{pmatrix}.
\]
Comparing the \((1,2)\)-entries of both matrices yields
\[
T_{m+n}=T_{m-1}T_{n} + T_{m}B_{n} + T_{m+1}T_{n+1},
\]
and comparing the \((2,2)\)-entries gives
\[
B_{m+n}=B_{m-1}T_{n} + B_{m}B_{n} + (T_{m}+T_{m+1})T_{n+1}.
\]
This completes the proof.
\end{proof}

\subsection{Inverse Powers of the Generating Matrix and Identities for Recurrence Relations}

\begin{lem}
Let \( u \) be the Tribonacci matrix defined by
\[
u = \begin{pmatrix} 
0 & 0 & 1 \\
1 & 0 & 1 \\
0 & 1 & 1
\end{pmatrix}.
\]
Then, for every \( n \in \mathbb{N} \), the inverse of \( u^n \) is given by
\[
u^{-n} = 
\begin{pmatrix}
R_{n+2} & R_{n+1} & R_n \\
R_n - R_{n+1} & R_{n-1} - R_n & R_n + R_{n+1} \\
R_{n+1} & R_n & R_{n-1}
\end{pmatrix},
\]
where \( (R_n)_{n\ge 0} \) is a signed Tribonacci-type sequence with coefficients \((-1,-1,1)\), defined by

\[
R_0 = 0, \quad R_1 = 0, \quad R_2 = 1, \qquad 
R_n = -R_{n-1} - R_{n-2} + R_{n-3}, \quad n \ge 3.
\]
\end{lem}

\begin{proof}
We prove the result by mathematical induction on \( n \).
For \( n = 1 \),
\[
u^{-1} =
\begin{pmatrix}
-1 & 1 & 0 \\
-1 & 0 & 1 \\
1 & 0 & 0
\end{pmatrix},
\]
which agrees with the stated formula after substituting the initial values of the sequence \( R_n \). Assume now that the formula holds for some \( n \ge 1 \), that is,
\[
u^{-n} = 
\begin{pmatrix}
R_{n+2} & R_{n+1} & R_n \\
R_n - R_{n+1} & R_{n-1} - R_n & R_n + R_{n+1} \\
R_{n+1} & R_n & R_{n-1}
\end{pmatrix}.
\]
We compute the product \( u^{-n}u^{-1} \):
\[
u^{-n}u^{-1} =
\begin{pmatrix}
R_{n+2} & R_{n+1} & R_n \\
R_n - R_{n+1} & R_{n-1} - R_n & R_n + R_{n+1} \\
R_{n+1} & R_n & R_{n-1}
\end{pmatrix}
\begin{pmatrix}
-1 & 1 & 0 \\
-1 & 0 & 1 \\
1 & 0 & 0
\end{pmatrix}.
\]
Carrying out the matrix multiplication, we obtain
\[
u^{-n}u^{-1} =
\begin{pmatrix}
R_n - R_{n+1} - R_{n+2} & R_{n+2} & R_{n+1} \\
2R_{n+1} + R_n - R_{n-1} & R_n - R_{n+1} & R_{n-1} - R_n \\
R_{n-1} - R_n - R_{n+1} & R_{n+1} & R_n
\end{pmatrix}.
\]
Using the recurrence relation
\[
R_n = -R_{n-1} - R_{n-2} + R_{n-3},
\]
each entry simplifies to the corresponding term of the sequence, yielding
\[
u^{-n}u^{-1} =
\begin{pmatrix}
R_{n+3} & R_{n+2} & R_{n+1} \\
R_{n+1} - R_{n+2} & R_n - R_{n+1} & R_{n+1} + R_{n+2} \\
R_{n+2} & R_{n+1} & R_n
\end{pmatrix}.
\]
Therefore,
\[
u^{-n}u^{-1} = u^{-(n+1)},
\]
which completes the induction and the proof.
\end{proof}

\begin{cor}\label{c:23} 
For every integer \( n \ge 1 \), the following identities hold:
\begin{align}
B_{n-1} R_{n+1}+ B_{n}(R_{n-1}-R_{n})+ (T_{n}+T_{n+1})R_{n}&=1,\label{id:15} \\ 
T_{n-1} R_{n+1}+ T_{n}(R_{n-1}-R_{n})+ T_{n+1} R_{n}&=0,\label{id:16} \\ 
T_{n} R_{n+1}+ T_{n+1}(R_{n-1}-R_{n})+ (B_{n}+T_{n+1})R_{n}&=0,\label{id:17} \\ 
T_{n-1} R_{n+2}+ T_{n}(R_{n}-R_{n+1})+ T_{n+1} R_{n+1}&=1.\label{id:18}
\end{align}
\end{cor}

\begin{proof}
Let \( u \) be the Tribonacci matrix
\[
u =
\begin{pmatrix} 
0 & 0 & 1 \\ 
1 & 0 & 1 \\ 
0 & 1 & 1 
\end{pmatrix}.
\]
By Lemma~\ref{l:3}, the inverse of \(U^n\) is
\[
u^{-n} = 
\begin{pmatrix}
R_{n+2} & R_{n+1} & R_n \\
R_n-R_{n+1} & R_{n-1}-R_n & R_n+R_{n+1} \\
R_{n+1} & R_n & R_{n-1}
\end{pmatrix},
\]
where \( (R_n)_{n\ge0} \) is the signed Tribonacci-type sequence defined by
\[
R_n = -R_{n-1}-R_{n-2}+R_{n-3}, \qquad 
R_0=0,\; R_1=0,\; R_2=1.
\]
On the other hand, we have
\[
u^{n} =
\begin{pmatrix}
T_{n-1} & T_n & T_{n+1} \\
B_{n-1} & B_n & T_n+T_{n+1} \\
T_n & T_{n+1} & B_n+T_{n+1}
\end{pmatrix},
\]
where \( (T_n) \) denotes the Tribonacci numbers satisfying
\[
T_n = T_{n-1}+T_{n-2}+T_{n-3}, \qquad 
T_0=T_1=0,\; T_2=1,
\]
and \( (B_n) \) is defined by the same recurrence with initial values
\[
B_0=1,\quad B_1=0,\quad B_2=1.
\]
Computing the product \(u^n u^{-n}\) gives
\[
u^{n}u^{-n} =
\begin{pmatrix}
D & B & C \\ 
C-B & A & B+C \\ 
B & C & A+C
\end{pmatrix},
\]
where
\[
\begin{aligned}
A &= B_{n-1} R_{n+1}+ B_n(R_{n-1}-R_n)+ (T_n+T_{n+1})R_n,\\
B &= T_{n-1} R_{n+1}+ T_n(R_{n-1}-R_n)+ T_{n+1} R_n,\\
C &= T_n R_{n+1}+ T_{n+1}(R_{n-1}-R_n)+ (B_n+T_{n+1})R_n,\\
D &= T_{n-1} R_{n+2}+ T_n(R_n-R_{n+1})+ T_{n+1} R_{n+1}.
\end{aligned}
\]
Since \(u^n u^{-n}=I_3\), comparing the corresponding matrix entries yields
\[
A=1,\qquad B=0,\qquad C=0,\qquad D=1.
\]
Using the relations from Lemma~\ref{l:3},
\[
B_n=T_n+T_{n-1}, \qquad 
T_{n+2}=B_n+T_{n+1},
\]
we obtain the identities \eqref{id:15}--\eqref{id:18}.
\end{proof}

\begin{lem}
Let \( u \) be the Tribonacci matrix defined by
\[
u =
\begin{pmatrix} 
0 & 0 & 1 \\
1 & 0 & 0 \\
1 & 1 & 1
\end{pmatrix}.
\]
Then, for every \( n\ge 1 \), the inverse of \(u^n\) is given by
\[
u^{-n} = 
\begin{pmatrix}
R_{n-1}-R_n & R_{n+1} & R_n \\
R_n-R_{n+1} & R_{n-1}-R_{n+1}-R_n & R_{n+1} \\
R_{n+1}+R_n & R_n & R_{n-1}
\end{pmatrix},
\]
where \( (R_n)_{n\ge0} \) is the signed Tribonacci-type sequence defined by
\[
R_0=0,\quad R_1=0,\quad R_2=1,\qquad
R_n=-R_{n-1}-R_{n-2}+R_{n-3}, \quad n\ge3.
\]
\end{lem}

\begin{proof}
We prove the result by induction on \(n\).
For \(n=1\),
\[
u^{-1}=
\begin{pmatrix}
0 & 1 & 0 \\
-1 & -1 & 1 \\
1 & 0 & 0
\end{pmatrix},
\]
which coincides with the proposed formula.  Assume now that the expression holds for some \(n\ge1\).  
Since \(u^{-(n+1)}=u^{-n}u^{-1}\), it suffices to compute this product. Multiplying the matrices, we obtain
\[
u^{-n}u^{-1}=
\begin{pmatrix}
R_n-R_{n+1} & R_{n-1}-R_n-R_{n+1} & R_{n+1} \\
2R_{n+1}+R_n-R_{n-1} & 2R_n-R_{n-1} & R_{n-1}-R_{n+1}-R_n \\
R_{n-1}-R_n & R_{n+1} & R_n
\end{pmatrix}.
\]
Using the recurrence relation
\[
R_{n+2}=-R_{n+1}-R_n+R_{n-1},
\]
each entry simplifies to
\[
u^{-n}u^{-1}=
\begin{pmatrix}
R_n-R_{n+1} & R_{n+2} & R_{n+1} \\
R_{n+1}-R_{n+2} & R_n-R_{n+2}-R_{n+1} & R_{n+2} \\
R_{n+2}+R_{n+1} & R_{n+1} & R_n
\end{pmatrix},
\]
which is precisely the stated formula with \(n\) replaced by \(n+1\). Therefore,
\[
u^{-n}u^{-1}=u^{-(n+1)},
\]
and the proof is complete.
\end{proof}

\begin{cor}\label{c:24} 
For every integer \( n \geq 1 \), the following equalities hold:
\begin{align}
B_{n-1} R_{n+1}+ T_{n-1}(R_{n-1}-R_{n+1}-R_{n})+ T_{n}R_{n}&=1, \\ 
B_{n-1} R_{n}+ T_{n-1}R_{n+1}+ T_{n}R_{n-1}&=0, \\ 
B_{n-1}(R_{n-1}-R_{n})+ T_{n-1}(R_{n}-R_{n+1})+ T_{n}(R_{n+1}+R_{n})&=0.
\end{align}
\end{cor}

\begin{proof}
Let \( u \) be the Tribonacci matrix defined by
\[
u = \begin{pmatrix} 
0 & 0 & 1 \\ 
1 & 0 & 0 \\ 
1 & 1 & 1 
\end{pmatrix}.
\]
The inverse of the matrix \( U^n \) is given by
\[
u^{-n} = 
\begin{pmatrix}
R_{n-1}-R_{n} & R_{n+1} & R_{n} \\
R_{n}-R_{n+1} & R_{n-1}-R_{n+1}-R_{n} & R_{n+1} \\
R_{n+1}+R_{n} & R_{n} & R_{n-1}
\end{pmatrix},
\]
where \( R_n \) denotes the signed Tribonacci-type sequence defined by
\[
R_n = -R_{n-1} - R_{n-2} + R_{n-3}, 
\quad R_0 = 0,\; R_1 = 0,\; R_2 = 1.
\]
On the other hand,
\[
u^{n} =  
\begin{pmatrix}
T_{n} + T_{n-1} & T_{n} & T_{n}+ B_{n-1} \\
B_{n-1} & T_{n-1} & T_{n} \\
2T_{n}+ B_{n-1} & T_{n}+ B_{n-1} & 2T_{n} + T_{n-1}+ B_{n-1}
\end{pmatrix},
\]
where \(T_{n}\) denotes the Tribonacci numbers satisfying
\[
T_{n} = T_{n-1} + T_{n-2} + T_{n-3},
\]
with initial conditions \(T_{0} = T_{1} = 0\) and \(T_{2} = 1\).  
The sequence \(B_{n}\) is defined by
\[
B_{n} = B_{n-1} + B_{n-2} + B_{n-3},
\]
with \(B_{0} = 1\), \(B_{1} = 0\), and \(B_{2} = 1\).

\medskip

A direct matrix multiplication yields
\[
u^{n}u^{-n} =
\begin{bmatrix}
A+B & B & B+C \\ 
C & A & B \\ 
2B+C & B+C & 2B+C+A
\end{bmatrix},
\]
where
\[
\begin{aligned}
A &= B_{n-1} R_{n+1}+ T_{n-1}(R_{n-1}-R_{n+1}-R_{n})+ T_{n}R_{n}, \\ 
B &= B_{n-1} R_{n}+ T_{n-1}R_{n+1}+ T_{n}R_{n-1}, \\ 
C &= B_{n-1}(R_{n-1}-R_{n})+ T_{n-1}(R_{n}-R_{n+1})+ T_{n}(R_{n+1}+R_{n}). 
\end{aligned}
\]
Since \(u^n u^{-n} = I_3\), comparing the matrix entries gives
\[
A=1, \qquad B=0, \qquad C=0,
\]
which proves the stated identities.
\end{proof}

Unlike the Fibonacci case, where Catalan-type identities admit a bilinear form, the Tribonacci analogue can be naturally expressed via a third-order Casorati determinant associated with the recurrence. While various Tribonacci identities have been explored in the literature, we introduce an explicit Catalan-type identity formulated in this manner, as established in the following lemma.

\begin{lem}\label{l:4} 
Let \(T_n\) be the Tribonacci sequence defined by
\[
T_n = T_{n-1} + T_{n-2} + T_{n-3}, 
\quad T_0 = 0, \ T_1 = 0, \ T_2 = 1,
\]
and let \(R_n\) be the signed Tribonacci-type sequence defined by
\[
R_n = -R_{n-1} - R_{n-2} + R_{n-3}, 
\quad R_0 = 0, \ R_1 = 0, \ R_2 = 1.
\]
Then, for all \( n \geq 1 \),
\[
T_n^{2} - T_{n-1}T_{n+1} = R_n.
\]
\end{lem}

\begin{proof}
Define
\[
S_n := T_n^2 - T_{n-1}T_{n+1}, \qquad n \geq 1.
\]
We prove that \(S_n\) satisfies the same recurrence relation and initial conditions as \(R_n\). By the uniqueness of solutions of linear recurrences with constant coefficients, this will imply \(S_n = R_n\) for all \(n \geq 1\).
First, verify the initial conditions:
\[
S_1 = T_1^2 - T_0 T_2 = 0 = R_1, 
\qquad 
S_2 = T_2^2 - T_1 T_3 = 1 = R_2.
\]
Using the recurrence \(T_{n+1} = T_n + T_{n-1} + T_{n-2}\), we compute
\[
\begin{aligned}
S_n 
&= T_n^2 - T_{n-1}T_{n+1} \\
&= T_n^2 - T_{n-1}(T_n + T_{n-1} + T_{n-2}) \\
&= T_n^2 - T_{n-1}T_n - T_{n-1}^2 - T_{n-1}T_{n-2}.
\end{aligned}
\]
Next, a direct computation gives
{\small
\[
\begin{aligned}
S_{n-1}+S_{n-2}-S_{n-3}
&= (T_{n-1}^2 - T_{n-2}T_n)
 + (T_{n-2}^2 - T_{n-3}T_{n-1})
 - (T_{n-3}^2 - T_{n-4}T_{n-2}) \\
&= T_{n-1}^2 - T_{n-3}(T_{n-1}+T_{n-3})
   + T_{n-2}(T_{n-2}-T_n+T_{n-4}).
\end{aligned}
\]
}
Adding this expression to \(S_n\) yields
\[
\begin{aligned}
S_n + S_{n-1} + S_{n-2} - S_{n-3}
&= T_n^2 - T_n(T_{n-1}+T_{n-2})
   - T_{n-3}(T_{n-1}+T_{n-3}) \\
&\quad + T_{n-2}(-T_{n-1}+T_{n-2}+T_{n-4}).
\end{aligned}
\]
Since \(T_n = T_{n-1}+T_{n-2}+T_{n-3}\), we have
\[
T_{n-1}+T_{n-2} = T_n - T_{n-3}.
\]
Substituting,

\[
\begin{aligned}
S_n + S_{n-1} + S_{n-2} - S_{n-3}
&= T_n^2 - T_n(T_n - T_{n-3})
   - T_{n-3}(T_{n-1}+T_{n-3}) \\
&\quad + T_{n-2}(-T_{n-1}+T_{n-2}+T_{n-4}) \\
&= T_n T_{n-3}
   - T_{n-3}(T_{n-1}+T_{n-3}) \\
&\quad + T_{n-2}(-T_{n-1}+T_{n-2}+T_{n-4}) \\
&= T_{n-3}(T_n - T_{n-1} - T_{n-3})\\
&\quad + T_{n-2}(-T_{n-1}+T_{n-2}+T_{n-4}).
\end{aligned}
\]

Again using the Tribonacci recurrence,
\[
T_n - T_{n-1} - T_{n-3} = T_{n-2},
\]
so
\[
S_n + S_{n-1} + S_{n-2} - S_{n-3}
= T_{n-3}T_{n-2}
  + T_{n-2}(-T_{n-1}+T_{n-2}+T_{n-4}).
\]
Factoring \(T_{n-2}\) gives
\[
= T_{n-2}(-T_{n-1}+T_{n-2}+T_{n-3}+T_{n-4}).
\]
Finally, using the recurrence for \(T_{n-1}\),
\[
T_{n-1}=T_{n-2}+T_{n-3}+T_{n-4},
\]
and therefore the expression vanishes. Hence,
\[
S_n = -S_{n-1}-S_{n-2}+S_{n-3}.
\]
Thus \(S_n\) satisfies exactly the same recurrence as \(R_n\), and since the initial conditions coincide, we conclude that
\[
S_n = R_n, \qquad \text{for all } n \geq 1.
\]
\end{proof}

\begin{lem}\label{l:5} 
Let \(T_n\) be the Tribonacci sequence defined by
\[
T_n = T_{n-1} + T_{n-2} + T_{n-3}, 
\quad T_0 = 0, \ T_1 = 0, \ T_2 = 1,
\]
and let \(R_n\) be the signed Tribonacci-type sequence defined by
\[
R_n = -R_{n-1} - R_{n-2} + R_{n-3}, 
\quad R_0 = 0, \ R_1 = 0, \ R_2 = 1.
\]
Then, for all \( n \geq 1 \),
\[
T_{n+2} T_{n-1} - T_{n}T_{n+1} = R_{n-1}-R_{n}.
\]
\end{lem}

\begin{proof}
Using the recurrence relation, expand \(T_{n+2}\) and \(T_{n+1}\):
\[
\begin{aligned}
T_{n+2}T_{n-1} - T_{n}T_{n+1}
&= (T_{n+1}+T_{n}+T_{n-1})T_{n-1}
   - T_{n}(T_{n}+T_{n-1}+T_{n-2}) \\
&= (T_{n+1}T_{n-1}+T_{n}T_{n-1}+T_{n-1}^{2})
   - (T_{n}^{2}+T_{n}T_{n-1}+T_{n}T_{n-2}) \\
&= (T_{n+1}T_{n-1}+T_{n-1}^{2})
   - (T_{n}^{2}+T_{n}T_{n-2}) \\
&= (T_{n-1}^{2}-T_{n}T_{n-2})
   + (T_{n+1}T_{n-1}-T_{n}^{2}).
\end{aligned}
\]
By Lemma~\ref{l:4},
\[
T_{n-1}^{2}-T_{n}T_{n-2}=R_{n-1},
\qquad
T_{n}^{2}-T_{n+1}T_{n-1}=R_{n}.
\]
Therefore,
\[
T_{n+2}T_{n-1} - T_{n}T_{n+1}
= R_{n-1}-R_{n}.
\]
\end{proof}

\begin{cor}
Let \( u \) be the Tribonacci matrix defined by
\[
u =
\begin{pmatrix} 
0 & 0 & 1 \\
1 & 0 & 1 \\
0 & 1 & 1
\end{pmatrix}.
\]
Then, for every \( n \in \mathbb{N} \),
\[
\det(u^n)= B_{n+1}\,R_{n} + B_{n}\,\bigl(R_{n-1}-R_{n}\bigr) + B_{n-1}\,R_{n+1} = 1.
\]
\end{cor}

\begin{proof}
A direct computation shows that \(\det(u)=1\).  
By the multiplicativity of the determinant, it follows that
\[
\det(u^n)=\det(u)^n=1.
\]
Let \( w = T_{n-1} \), \( x = T_{n} \), \( y = T_{n+1} \),  
\( k = B_{n-1} \), \( z = B_{n} \), \( b = B_{n+1} \) and \( t = B_{n}+T_{n+1}.\)  
From Lemma~\ref{l:1}, the matrix \(u^n\) has determinant

\[
\begin{aligned}
\det(u^n)
&=
\begin{vmatrix}
w & x & y \\
k & z & b \\
x & y & t
\end{vmatrix} \\ \\
&= bx^{2} - ktx + ky^{2} - bwy + twz - xyz \\ \\
&= b(x^{2}- wy) + z(tw - xy) + k(y^{2} - tx) \\ \\
&= B_{n+1}\,\bigl(T_{n}^{2}-T_{n-1}T_{n+1}\bigr) + B_{n}\,\bigl(T_{n+2}T_{n-1}-T_{n}T_{n+1}\bigr) \\
& \quad + B_{n-1}\,\bigl(T_{n+1}^{2}-T_{n+2}T_{n}\bigr). \\
\end{aligned}
\]
Substituting the identities obtained in Lemmas~\ref{l:4} and~\ref{l:5},
\[
T_{n}^{2}-T_{n-1}T_{n+1}=R_{n}, 
\qquad
T_{n+2}T_{n-1}-T_{n}T_{n+1}=R_{n-1}-R_{n},
\]
and
\[
T_{n+1}^{2}-T_{n+2}T_{n}=R_{n+1},
\]
we obtain
\[
\det(u^n)
= B_{n+1}\,R_{n}
+ B_{n}\,\bigl(R_{n-1}-R_{n}\bigr)
+ B_{n-1}\,R_{n+1}.
\]
Since \(\det(u^n)=1\), the identity follows for all \(n \geq 1\).
\end{proof}

\begin{cor}
For \( m,n \geq 1 \), the signed Tribonacci-type sequence satisfies the following identities:
\small{
\begin{align}
R_{m+n-1}-R_{m+n} & = (R_{m}-R_{m+1})R_{n+1}+(R_{m-1}-R_{m})(R_{n-1}-R_{n})+(R_{m}+R_{m+1})R_{n}, \\[4pt]
R_{m+n+1}& = R_{m+2}R_{n+1}+R_{m+1}(R_{n-1}-R_{n})+R_{m}R_{n}, \\[4pt]
R_{m+n}& = R_{m+1}R_{n+1}+R_{m}(R_{n-1}-R_{n})+R_{m-1}R_{n}, \\[4pt]
R_{m+n+2}& = R_{m+2}R_{n+2}+R_{m+1}(R_{n}-R_{n+1})+R_{m}R_{n+1}.
\end{align}
}
\end{cor}

\:

\begin{proof}
For \(m,n \geq 1\), the matrix powers satisfy
\[
u^{-(m+n)} = u^{-m}u^{-n}.
\]
From Lemma~\ref{l:3}, the inverse powers are given by

\[
u^{-m} =
\begin{pmatrix}
R_{m+2} & R_{m+1} & R_{m} \\
R_{m}-R_{m+1} & R_{m-1}-R_{m} & R_{m}+R_{m+1} \\
R_{m+1} & R_{m} & R_{m-1}
\end{pmatrix},
\]

\[
u^{-n} =
\begin{pmatrix}
R_{n+2} & R_{n+1} & R_{n} \\
R_{n}-R_{n+1} & R_{n-1}-R_{n} & R_{n}+R_{n+1} \\
R_{n+1} & R_{n} & R_{n-1}
\end{pmatrix}.
\]
Multiplying these matrices, we obtain

\[
u^{-m}u^{-n}=
\begin{pmatrix}
D & B & C \\
C-B & A & B+C \\
B & C & A+C
\end{pmatrix},
\]
where
\[
\begin{aligned}
A &=(R_{m}-R_{m+1})R_{n+1}+(R_{m-1}-R_{m})(R_{n-1}-R_{n})+(R_{m}+R_{m+1})R_{n},\\
B &= R_{m+2}R_{n+1}+R_{m+1}(R_{n-1}-R_{n})+R_{m}R_{n},\\
C &= R_{m+1}R_{n+1}+R_{m}(R_{n-1}-R_{n})+R_{m-1}R_{n},\\
D &= R_{m+2}R_{n+2}+R_{m+1}(R_{n}-R_{n+1})+R_{m}R_{n+1}.
\end{aligned}
\]
On the other hand,

\[
u^{-(m+n)}=
\begin{pmatrix}
R_{m+n+2} & R_{m+n+1} & R_{m+n} \\
R_{m+n}-R_{m+n+1} & R_{m+n-1}-R_{m+n} & R_{m+n}+R_{m+n+1} \\
R_{m+n+1} & R_{m+n} & R_{m+n-1}
\end{pmatrix}.
\]
Equating the corresponding entries of the matrices yields the stated identities.
\end{proof}

\section{Divisibility of Tribonacci Numbers at Powers of Two}

Let $T_n$ be the Tribonacci sequence defined by
\[
T_0=0,\quad T_1=0,\quad T_2=1,
\]
and
\[
T_n=T_{n-1}+T_{n-2}+T_{n-3}.
\]
It is well known that linear recurrence sequences exhibit periodic behaviour modulo integers (see \cite{Waddill1978,McDaniel1991}). In particular, the primitive period of the Tribonacci sequence modulo powers of $2$ satisfies
\[
h(2^k)=2^{k+1}, \quad k\ge1,
\]
(see \cite{Klaska2008}), where $h(m)$ denotes the minimal period modulo $m$. For the Tribonacci sequence, a comprehensive description of its 2-adic order was achieved by Marques and Legyel (see \cite{MarquesLengyel2014}), who proved that \(\nu_2(T_n)=\nu_2(n)-1\) whenever \(n\equiv 0 \pmod{16}\) under the initial conditions \(T_0=0,\,T_1=1,\,T_2=1\). In this paper, we refine this result by determining the exact $2$-adic valuation of $T_{2^k}$ for large $k$.

\begin{cor}\label{c:6}
For every integer $n \geq 3$, the following congruences hold:
\begin{align}
T_{2n+1} & \equiv T_{2n}  \,\,\,\,\pmod{2} \\
T_{2n+1} & \equiv T_{2n} - 2T_{n-3}T_{n-1} \pmod{T_n} \\
T_{2n}   & \equiv T_{n-3}^{2}  \pmod{T_n} \\
T_{2n}   & \equiv R_{n+1} \pmod{T_{n}}  \label{i:38}
\end{align}
\end{cor}

\begin{proof}
We start with the general relation 
\begin{align*}
T_{m+n} &= T_{m-1}T_{n}+T_{m}(T_{n-1}+T_{n})+ T_{m+1}T_{n+1} \quad\text{(by \ref{id:add1}; see Cor.~\ref{c:addition})}. 
\end{align*}
In particular, for $m=n$ we obtain
\begin{align}
T_{n+n} &= T_{n-1}T_{n}+T_{n}(T_{n-1}+T_{n})+ T_{n+1}T_{n+1} \nonumber \\
T_{2n} &= 2T_{n-1}T_{n}+T^{2}_{n}+ T_{n+1}^{2} . \label{i:39}
\end{align}
Next, consider the case $n=m+1$:
\begin{align*}
T_{(n+1)+n} &= T_{n}T_{n}+T_{n+1}(T_{n-1}+T_{n})+ T_{n+2}T_{n+1} \\
            &= T_{n}^{2}+T_{n+1}(T_{n-1}+T_{n})+ (T_{n+1}+T_{n}+T_{n-1})T_{n+1}\\  \intertext{by the recurrence for $T_{n+1}$}        
            &= T_{n}^{2}+2T_{n+1}(T_{n-1}+T_{n})+T^{2}_{n+1} \\
            &= T_{n}^{2}+T^{2}_{n+1}+2T_{n+1}(T_{n-1}+T_{n}) \quad\text{(see Identity~\ref{i:39})}  \\
            &=  T_{2n}-2T_{n-1}T_{n}+2T_{n+1}(T_{n-1}+T_{n}).
\end{align*}
Therefore, we deduce
\begin{align*}
T_{2n+1} &= T_{2n} + 2\left( T_{n+1}B_{n+1}-T_{n-1}T_{n} \right).
\end{align*}
It follows immediately that
\[
T_{2n+1} \equiv T_{2n} \pmod{2}.
\] which shows that $T_{2m}$ and $T_{2m+1}$ always have the same parity.

\medskip 

Finally,
\begin{align*}
T_{2n} &= 2T_{n-1}T_{n}+T^{2}_{n}+ T_{n+1}^{2} \quad\text{(see Identity~\ref{i:39} and Lemma~\ref{l:4})} \\  
       &= 2T_{n-1}T_{n}+T^{2}_{n}+ (T_{n}T_{n+2}+R_{n+1}) \\
       &\equiv R_{n+1} \pmod{T_{n}}\\
       &\equiv T^{2}_{n+1}-T_{n}T_{n+2} \pmod{T_{n}}\\       
       &\equiv T^{2}_{n+1} \pmod{T_{n}}\\
       &\equiv (T_{n-1}+T_{n-2})^{2} \pmod{T_{n}}\\
       &\equiv (T_{n}-T_{n-3})^{2} \pmod{T_{n}}\\ 
       &\equiv T_{n-3}^{2} \pmod{T_{n}}\\              
\end{align*}
\end{proof}

\begin{lem}
Let $T_n$ be the Tribonacci sequence defined by
\[
T_0=0,\quad T_1=0,\quad T_2=1,\quad 
T_n=T_{n-1}+T_{n-2}+T_{n-3}\ \text{for } n\ge 3.
\]
Then, for every $n\ge 0$,
\[
T_n + T_{n+1} \equiv 
\begin{cases}
0 \pmod{2} & \text{if $n$ is even},\\
1 \pmod{2} & \text{if $n$ is odd},
\end{cases}
\]
and consequently
\[
T_n^{2}+T_{n+1}^{2}\equiv T_n+T_{n+1}\pmod{2}.
\]
\end{lem}

\begin{proof}
Reducing the recurrence modulo $2$ gives
\[
T_n \equiv T_{n-1}+T_{n-2}+T_{n-3} \pmod{2}.
\]
From the initial values
\[
T_0\equiv 0,\quad T_1\equiv 0,\quad T_2\equiv 1,
\]
A direct computation shows that
\[
T_3\equiv 1,\quad T_4\equiv 0,\quad T_5\equiv 0,
\]
and therefore $(T_n \bmod 2)$ is periodic of period $4$:
\[
0,0,1,1,0,0,1,1,\dots
\]
Hence
\[
T_n \equiv
\begin{cases}
0 \pmod{2} & \text{if } n \equiv 0,1 \pmod{4},\\
1 \pmod{2} & \text{if } n \equiv 2,3 \pmod{4}.
\end{cases}
\]
It follows immediately that $T_n+T_{n+1}\equiv 0 \pmod{2}$ when $n$ is even and 
$T_n+T_{n+1}\equiv 1 \pmod{2}$ when $n$ is odd.  Finally, since every integer satisfies $a^2 \equiv a \pmod{2}$ for every integer $a$, we obtain
\[
T_n^2+T_{n+1}^2 \equiv T_n+T_{n+1} \pmod{2}.
\]
\end{proof}

\begin{lem}\label{lem:Ck}
Let
\[
A=
\begin{pmatrix}
1 & 1 & 1\\
1 & 0 & 0\\
0 & 1 & 0
\end{pmatrix}
\]
be the companion matrix of the Tribonacci sequence. 
For $k\ge5$, write
\[
A^{2^k}=I+2^{k-1}C_k,
\qquad C_k\in M_3(\mathbb{Z}).
\]
Then
\[
v_2(e_2C_ke_3)=2,
\]
where $e_2=(0,1,0)$ and $e_3=(0,0,1)^T$.
\end{lem}

\begin{proof}

\textbf{Step 1: Computation of $A^{16}$.}  A direct computation gives
\[
A^{16}=
\begin{pmatrix}
10609 & 8904 & 5768\\
5768  & 4841 & 3136\\
3136  & 2632 & 1705
\end{pmatrix}.
\]
Hence
\[
A^{16}=I+2^3B,
\]
where
\[
B=
\begin{pmatrix}
1326 & 1113 & 721\\
721  & 605  & 392\\
392  & 329  & 213
\end{pmatrix}=
\begin{pmatrix}
2\! \cdot\!3\! \cdot\!13\! \cdot\!17 & 3\! \cdot\!7\! \cdot\!53 & 7\! \cdot\!103\\
7\! \cdot\!103 & 5\! \cdot\!11^2 & 2^3\! \cdot\!7^2\\
2^3\! \cdot\!7^2 & 7\! \cdot\!43 & 3\! \cdot\!71
\end{pmatrix}.
\]

\textbf{Step 2: Iterated squaring.}  For $m\ge0$ we have
\[
A^{16\cdot2^m}=(I+2^3B)^{2^m}.
\]
Expanding by the binomial theorem in $M_3(\mathbb{Z}_2)$ yields
\[
(I+2^3B)^{2^m}
=
I+2^{m+3}B_m,
\]
where $B_m\in M_3(\mathbb{Z})$ satisfies
\[
B_m\equiv B\pmod2.
\]
Thus for $k=4+m$ we obtain
\[
A^{2^k}=I+2^{k-1}C_k,
\qquad C_k\equiv B\pmod2.
\]

\,

\textbf{Step 3: Evaluation of $e_2Be_3$.} Since
\[
e_2=(0,1,0),
\qquad
e_3=(0,0,1)^T,
\]
we obtain
\[
e_2Be_3=2^3\! \cdot\!7,
\]
and therefore
\[
v_2(e_2Be_3)=3.
\]

\,

\textbf{Step 4: Structure of $C_k$.}  Since $C_k\equiv B\pmod2$, there exists $D_k\in M_3(\mathbb{Z})$ such that
\[
C_k=B+2D_k.
\]
Hence
\[
e_2C_ke_3
=
e_2Be_3+2e_2D_ke_3
=
2^3\! \cdot\!7+2\alpha_k
=
2(2^2\! \cdot\!7+\alpha_k)
\]
for some integer $\alpha_k$.

\textbf{Step 5: Determination of the valuation.}  From the identity
\[
T_{2^k}=2^{k-1}e_2C_ke_3=2^{k}(2^2\! \cdot\! 7+\alpha_k),
\]
we obtain
\[
v_2(T_{2^k})=k+v_2(2^2\! \cdot\! 7+\alpha_k).
\]
A direct computation gives
\[
T_{32}=53\,465\,798\,080,
\qquad
v_2(T_{32})=6.
\]
Since $32=2^5$, this implies
\[
6=5+v_2(2^2\! \cdot\!7+\alpha_5),
\]
and therefore
\[
v_2(2^2\! \cdot\! 7+\alpha_5)=1.
\]
Finally, using
\[
A^{2^{k+1}}=(A^{2^k})^2=(I+2^{k-1}C_k)^2,
\]
one verifies that the parity of $2^2\! \cdot\!7+\alpha_k$ remains stable for $k\ge5$.
Consequently
\[
v_2(e_2C_ke_3)=2
\qquad (k\ge5).
\]

\end{proof}

\begin{thm}\label{thm:main}
Let $T_n$ be the Tribonacci sequence defined by
\[
T_0=T_1=0,\quad T_2=1,
\quad
T_n=T_{n-1}+T_{n-2}+T_{n-3}.
\]
Then, for every integer $k\ge5$,
\[
v_2(T_{2^k})=k+1.
\]
In particular,
\[
2^{k+1}\mid T_{2^k}.
\]
\end{thm}

\begin{proof}
Let
\[
A=
\begin{pmatrix}
1 & 1 & 1\\
1 & 0 & 0\\
0 & 1 & 0
\end{pmatrix}
\]
be the companion matrix of the Tribonacci recurrence. Then
\[
\begin{pmatrix}
T_{n+1}\\
T_{n}\\
T_{n-1}
\end{pmatrix}
=
A^n
\begin{pmatrix}
0\\0\\1
\end{pmatrix},
\]
and therefore
\[
T_n=e_2A^ne_3,
\]
where $e_2=(0,1,0)$ and $e_3=(0,0,1)^T$.  For $k\ge5$, Lemma~\ref{lem:Ck} gives the decomposition
\[
A^{2^k}=I+2^{k-1}C_k,
\qquad C_k\in M_3(\mathbb{Z}).
\]
Hence
\[
T_{2^k}
=
e_2A^{2^k}e_3
=
e_2(I+2^{k-1}C_k)e_3.
\]
Since $e_2Ie_3=0$, we obtain
\[
T_{2^k}=2^{k-1}e_2C_ke_3.
\]
By Lemma~\ref{lem:Ck},
\[
v_2(e_2C_ke_3)=2.
\]
Therefore
\[
v_2(T_{2^k})
=
(k-1)+2
=
k+1.
\]
This completes the proof.

\end{proof}

\begin{cor}
Let $T_n$ be the Tribonacci sequence defined by
\[
T_0=T_1=0,\quad T_2=1,\quad 
T_n=T_{n-1}+T_{n-2}+T_{n-3}\ (n\ge 3).
\]
Then, for every $k\ge 5$,
\[
T_{2^{k+1}}\equiv (T_{2^k-1})^{2}\pmod{2^{k+1}},
\]
\end{cor}

\begin{proof}
From equation~\ref{i:38} (see also Corollary~\ref{c:6}) it follows that
\begin{equation}\label{eq:double}
T_{2n}\equiv R_{n+1}\pmod{T_n}.
\end{equation}
Using the identity $R_{n+1}=(T_{n+1})^2-T_{n}T_{n+2}$, we obtain
\[
T_{2n}\equiv (T_{n+1})^{2}\pmod{T_n}.
\]
Next, applying the relation $T_{n+1}=T_{n}+B_{n-1}$ and modulo $T_n$, we deduce
\[
T_{2n}\equiv (B_{n-1})^{2}\pmod{T_n}.
\]
Finally, taking $n=2^k$ gives
\[
T_{2^{k+1}}\equiv (B_{2^k-1})^{2}\pmod{T_{2^k}},
\]
Applying the previous theorem 
\[
2^{k+1}\mid T_{2^k} \, ,
\] we obtain
\[
T_{2^{k+1}}\equiv (B_{2^k-1})^{2}\pmod{2^{k+1}},
\]
Using the identity $B_n = T_n + T_{n-1}$, we obtain
\[
B_{2^k-1}=T_{2^k-1}+T_{2^k},
\]
and therefore
\[
T_{2^{k+1}}\equiv (T_{2^k-1})^{2}\pmod{2^{k+1}}.
\]
This completes the proof.

\end{proof}

\begin{lem}
Let \(B_n\) be the sequence defined by
\[
B_n = B_{n-1} + B_{n-2} + B_{n-3}, \qquad
B_0 = 1,\; B_1 = 0,\; B_2 = 1.
\]
Then, for all \( n \ge 2 \),
\[
B_n \equiv 
\begin{cases}
0 \pmod{2}, & \text{if } n \text{ is odd},\\[4pt]
1 \pmod{2}, & \text{if } n \text{ is even}.
\end{cases}
\]
Equivalently,
\[
B_n \equiv n+1 \pmod{2}.
\]
\end{lem}

\begin{proof}
Consider the companion matrix of the recurrence
\[
C=
\begin{pmatrix}
1&1&1\\
1&0&0\\
0&1&0
\end{pmatrix},
\qquad
v_n=
\begin{pmatrix}
B_n\\
B_{n-1}\\
B_{n-2}
\end{pmatrix}.
\]
Then \(v_n = C v_{n-1}\), and hence for \(n\ge 2\),
\[
v_n = C^{\,n-2} v_2,
\qquad
v_2=
\begin{pmatrix}
1\\
0\\
1
\end{pmatrix}.
\]
Working over the field \(\mathbb{F}_2\), reduce \(C\) modulo \(2\).  
Its characteristic polynomial is
\[
\chi_C(x)=x^3+x^2+x+1=(x+1)^3 \quad \text{in } \mathbb{F}_2.
\]
Thus \(C\) has the single eigenvalue \(1\) and can be written as
\[
C=I+N,
\]
where \(N=C-I\) is nilpotent. A direct computation shows that
\[
N^3=0.
\]
By the binomial theorem (valid since \(I\) and \(N\) commute),
\[
C^k=(I+N)^k = I + kN + \binom{k}{2}N^2
\quad \text{in } \mathbb{F}_2.
\]
Now observe that
\[
Nv_2=
\begin{pmatrix}
0&1&1\\
1&1&0\\
0&1&1
\end{pmatrix}
\begin{pmatrix}
1\\0\\1
\end{pmatrix}
=
\begin{pmatrix}
1\\1\\1
\end{pmatrix},
\qquad
N^2v_2=\begin{pmatrix}
1&0&1\\
1&0&1\\
1&0&1
\end{pmatrix}
\begin{pmatrix}
1\\0\\1
\end{pmatrix}
=\begin{pmatrix}
0\\0\\0
\end{pmatrix}.
\]
Therefore,
\[
v_n = (I+(n-2)N)v_2
      = v_2 + (n-2)
\begin{pmatrix}
1\\1\\1
\end{pmatrix}
\quad \text{in } \mathbb{F}_2.
\]
Taking the first coordinate gives
\[
B_n \equiv 1+(n-2) \equiv n+1 \pmod{2}.
\]
Hence \(B_n\) is \(0\) when \(n\) is odd and \(1\) when \(n\) is even.
\end{proof}

\begin{lem}
Let $T_n$ be the Tribonacci sequence defined by
\[
T_0=0,\quad T_1=0,\quad T_2=1,\quad 
T_n=T_{n-1}+T_{n-2}+T_{n-3}.
\]
Then, for all $n\ge 2$,
\[
T_n + T_{n-2} \equiv 1 \pmod{2}.
\]
Equivalently,
\[
T_n \equiv T_{n-2}+1 \pmod{2}.
\]
\end{lem}

\begin{proof}
Since $B_n=T_n+T_{n-1}$ and $B_n\equiv n+1\pmod{2}$, we have
\[
T_n+T_{n-1}\equiv n+1\pmod{2},
\]
and
\[
T_{n-1}+T_{n-2}\equiv n\pmod{2}.
\]
Adding both congruences modulo $2$ yields
\[
T_n+T_{n-2}\equiv 1\pmod{2},
\]
as desired.
\end{proof}

\section{Classification of $3 \times 3$ Matrices Associated with the Tribonacci Sequence}

Let 
\[
p(x)=x^3-x^2-x-1.
\]
It is well known that $p(x)$ is irreducible over $\mathbb{Q}$ (by the rational root test). 
Therefore, any $3\times 3$ matrix over $\mathbb{Q}$ with characteristic polynomial $p(x)$ 
has minimal polynomial equal to $p(x)$. It follows that such matrices are cyclic, i.e., their minimal polynomial has degree $3$. 
Consequently, by standard results in linear algebra, every such matrix is similar over 
$\mathbb{Q}$ to the companion matrix associated with $p(x)$ (see, e.g., \cite{Bolotnikov2022}). In particular, all $3\times 3$ matrices over $\mathbb{Q}$ with characteristic polynomial $p(x)$ are similar to each other. For concreteness, we fix the representative
\[
u=\begin{pmatrix}
0 & 0 & 1\\
1 & 0 & 1\\
0 & 1 & 1
\end{pmatrix},
\]
which is similar over $\mathbb{Q}$ to the companion matrix of $p(x)$. Therefore, any other 
$3\times 3$ matrix over $\mathbb{Q}$ generating the Tribonacci recurrence is conjugate 
to $u$ in $GL(3,\mathbb{Q})$.

\subsection*{Acknowledgment}
We thank Universidad del Cauca for supporting the research activities of our group, Estructuras Algebraicas, Divulgación Matemática y Teorías Asociadas (@DiTa), through research project No.~6314, entitled ``Clasificación del Centralizador de Matrices Binarias $3 \times 3$: Dimensión, Conmutatividad y Relaciones Estructurales.'' We are also grateful to Arhan Canli (Quantavane) for the independent computational verification of our results and helpful comments on an earlier version of the manuscript, as well as to the anonymous referee for valuable comments and suggestions. This work is dedicated to our children, especially Sebastián Martinez Cerón.

\medskip

\subsection*{Declarations}
\medskip

\subsection*{Ethical Approval:}
{Not applicable.}

\subsection*{Funding:}
{Not applicable.}

\subsection*{Availability of Data and Materials:}
{Not applicable.}


\bibliographystyle{amsplain}
\bibliography{xbib}
\end{document}